\documentclass[12pt]{article}
\usepackage{amsfonts,amsmath,amssymb,mathrsfs}

\title{Global Existence of Bell's Time-Inhomogeneous Jump Process
for Lattice Quantum Field Theory}
\author{
    Hans-Otto Georgii\footnote{Mathematisches Institut der Universit\"{a}t
          M\"{u}nchen, Theresienstra{\ss}e 39, 80333 M\"{u}nchen, Germany.
          E-mail: georgii@mathematik.uni-muenchen.de}\ \ and
    Roderich Tumulka\footnote{Current address: Dipartimento di Fisica,
          Universit\`a di Genova, Via Dodecaneso 33, 16146 Genova, Italy.
          E-mail: tumulka@mathematik.uni-muenchen.de}
}
\date{March 31, 2004}

\newcommand{\CCC}{\mathbb{C}} 
\newcommand{\RRR}{\mathbb{R}} 
\newcommand{\R}{\mathbb{R}} 
\newcommand{\NNN}{\mathbb{N}} 
\newcommand{\ZZZ}{\mathbb{Z}} 
\newcommand{\EEE}{\mathbb{E}} 
\newcommand{\PPP}{\mathbb{P}} 
\newcommand{\E}{\mathrm{e}} 
\newcommand{\I}{\mathrm{i}} 
\newcommand{\1}{1}
\newcommand{\tr}{\mathrm{tr}} 
\renewcommand{\Im}{\mathrm{Im}} 
\newcommand{\rmP}{\mathrm{P}}
\newcommand{\vr}{{\boldsymbol r}}
\newcommand{\Hilbert}{\mathscr{H}}
\renewcommand{\sp}[2]{\langle #1 | #2 \rangle}
\newcommand{\conf}{E} 
\newcommand{\prob}{\PPP}
\newcommand{\pov}{{P}}
\newcommand{\D}{\mathscr{D}}
\newcommand{\Indexset}{\mathscr{I}}

\newcommand{\measure}{\pi} 
\newcommand{\friedhof}{\vartriangle}
\newcommand{\dest}{Z}
\newcommand{\node}{\mathcal{N}}

\newcommand{\x}{x}\newcommand{\y}{y}
\newcommand{\g}{\gamma}
\newcommand{\G}{\Gamma}
\renewcommand{\tau}{\theta}

\newtheorem{theorem}{Theorem}
\newtheorem{proposition}{Proposition}
\newtheorem{lemma}{Lemma}

\newtheorem{ass}{Assumption}

\newcommand{\Proof}[1][]{\noindent{\itshape Proof#1: }}
\newcommand{\EndProof}{~$\hfill\square$\bigskip}

\begin{document}\maketitle
\begin{abstract}
   We consider the time-inhomogeneous Markovian jump process introduced by
   John~S.~Bell~\cite{BellBeables} for a lattice quantum field theory,
   which runs on the associated configuration space.
   Its jump rates, tailored to give the process the quantum distribution
   $|\Psi_t|^2$ at all times $t$, typically exhibit singularities.
   We establish the existence of a unique such
   process for all times, under suitable assumptions on the
   Hamiltonian or the initial state vector $\Psi_0$.  The proof of
   non-explosion takes advantage of the special role of the $|\Psi_t|^2$
   distribution.

\medskip

\noindent
{Key words. } Markov jump processes, non-explosion, time-dependent jump rates,
 equivariant distributions, Bell's process, lattice quantum field theory.

\noindent
{MSC (2000). } \underline{60J75} (jump processes),
81T25 (quantum field theory on lattices)

\noindent
{PACS. } 03.65.Ta (foundations of quantum mechanics),
02.50.-r (probability theory, stochastic processes, and statistics),
03.70.+k (theory of quantized fields)
%
\end{abstract}

\section{Introduction}

This paper deals with the existence of Markov jump processes on 
countable sets having time-inhomogeneous transition rates of a 
particular
form  proposed by John S.~Bell \cite{BellBeables}
in his observer-independent formulation of lattice quantum field theories.
Bell's setting is as follows.

\bigskip \noindent\emph{Bell's model. } Let $\conf$ be the configuration
space for a variable but finite number of particles on a countable lattice
$\Lambda\subseteq \RRR^3$. A configuration $\x \in
\conf$ is mathematically represented by a function $\x:\Lambda \to
\ZZZ_+:=\{ 0,1,2, \ldots\}$ indicating the number of particles
$\x(\vr)$ at a site $\vr \in \Lambda$.
Thus
\[
    \conf = \Bigl\{\x\in \ZZZ_+^\Lambda: \textstyle\sum\limits_{\vr \in
    \Lambda} \x(\vr ) <\infty \Bigr\}.
\]
Hence $\conf$ is countably infinite. (In
Bell's proposal, $\x(\vr)$ is the number of \emph{fermions} at site $\vr$,
but this is of no relevance here.)

Next, Bell considers the Hilbert space $\Hilbert$ and the Hamiltonian $H$
of a lattice quantum field theory. This means that $H$ is a 
self-adjoint operator on $\Hilbert$ determining the quantum state
at time $t$ via
\begin{equation}\label{Schroed}
    \Psi_t = \E^{-\I H t/\hbar} \, \Psi_0
\end{equation}
for some initial state vector $\Psi_0$. These quantities are related
to the configuration space by a projection-valued measure (PVM) $\pov$
on $\conf$ acting on $\Hilbert$. That is, for every $\x\in \conf$
there exists an associated projection $\pov(\x)$ such that $\sum_{\x
  \in \conf} \pov(\x) = I$,  where $I$ is the identity operator, and
$\pov(\x)\pov(\y) =0$ when $\x \neq \y$. Specifically, $\pov(\x)$ is
the projection to the joint eigenspace of the (commuting)  fermion number
operators $N(\vr)$ associated with the eigenvalues $\x(\vr)$.  In
particular, $\sp{\Psi_t}{\pov(\x) \Psi_t}$ is the quantum probability
of a configuration $x$ at time $t$.
%
%
Bell then introduces the transition rate
\begin{equation}\label{tranrates}
    \sigma_t(\y|\x) = \frac{[(2/\hbar) \, \Im \, \sp{\Psi_t}{\pov(\y) H
    \pov(\x) \Psi_t}]^+}{\sp{\Psi_t}{\pov(\x) \Psi_t}}
\end{equation}
for a jump from $\x$ to $\y\in \conf$, where $a^+ = \max(a,0)$
denotes the positive part of $a$.
Note that $\sigma_t(\x|\x)=0$ because $\sp{\Psi_t}{\pov(\x) H \pov(\x)
\Psi_t} = \sp{\pov(\x) \Psi_t} { H \pov(\x) \Psi_t}$ is real. A formal
calculation yields  that
this choice of the jump rates is compatible with the process having
distribution $\sp{\Psi_t}{\pov(\,\cdot\,) \Psi_t}$
at each time $t$. See \cite{crea2} for an extensive discussion
of this jump rate formula.
In this paper we will choose the time unit
such that $\hbar=2$.

\bigskip\noindent\emph{Probabilistic questions. }
One of the main features of the transition rates \eqref{tranrates} is
that they become singular at times $t$ when $x$ becomes a ``node'' of
$\Psi_t$\,, i.e., when the denominator $\sp{\Psi_t} {\pov(\x) \Psi_t}$
in \eqref{tranrates} vanishes. So, at such times the process would not
know how to proceed. Fortunately, it turns out that the
increase of the rates close to such singularities has the positive effect
of forcing the process to jump away
before the singularity time is reached.

A more serious problem is the possibility of explosion
in finite time; that is, the jump times $T_n$ could accumulate so that
$\zeta = \sup_n T_n<\infty$ with positive probability.
The standard criteria
for non-explosion of pure jump processes are confined to
transition rates that are homogeneous in time, relying heavily
on the fact that the holding times are then exponentially distributed 
and independent; see, e.g., Section~2.7 of \cite{norris} or 
Proposition 10.21 of \cite{Kall}.  This
independence, however, fails to hold in the case of time-dependent 
jump rates, and the singularities  of Bell's transition rates
do not allow any simple bounds excluding explosion. The only thing one
knows is that the process is designed to have the prescribed quantum 
distribution at
fixed (deterministic) times, and it is this fact we will exploit.

Our proof will not make any
use of the particular construction or meaning of $\conf$ and $\pov$.
We will merely assume that $\conf$ is a countable set and $\pov$ a PVM
on $\conf$ acting on $\Hilbert$.  Actually we only need that
$\pov$ is a positive-operator-valued measure; see Section
\ref{sec:Result} below. Steps towards
an existence proof for Bell's process have already been made by
Bacciagaluppi \cite{Bthesis,BD}; his approach is, however, very
different from ours.

\bigskip\noindent\emph{Physical Perspective. }
Bell's observer-independent formulation of lattice quantum field theories
has attracted increasing attention recently
\cite{Colin,Dennis2,Dennis,crlet,crea1,crea2}.
Apart from its
relevance to the foundations of quantum theory, it has proven useful
for numerical simulations \cite{Dennis2}, and has been found
distinguished among all $|\Psi|^2$ distributed processes as the
minimal one \cite{crea2,Roy}, involving the least amount of
stochasticity.

There are close connections between Bell's model and two well-known
$|\Psi|^2$ distributed processes associated with nonrelativistic
quantum mechanics in $\RRR^3$: E.~Nelson's \emph{stochastic mechanics}
\cite{stoch1,stoch2,carlen1,carlen2} and \emph{Bohmian mechanics}
\cite{Bohm52,Bellhidden,survey,bmex}.  These processes are similar in
spirit to Bell's process, and can be combined with Bell's stochastic
jumps to include particle creation and annihilation
\cite{crlet,crea1,crea2}. Bell's process has also been utilized for
modal interpretations of quantum theory \cite{BD}.  Bohmian mechanics
arises as the continuum limit of Bell's process for a suitable choice
of $H$ and $\conf$ \cite{Sudbery,Vink}, and in general the continuum
limit presumably resembles the combined Bell--Bohm model of
\cite{crea1}.  A generalization of Bell's jump rate \eqref{tranrates}
to continuum spaces $\conf$ is given in \cite{crea2}.  The global
existence problem of stochastic mechanics has been solved in
\cite{carlen1} (see also \cite{carlen2,stoch1}) and the one of Bohmian
mechanics in \cite{bmex}, whereas for combined models with jumps, such
as the ones considered in \cite{crea1,crea2}, it is still open. The
existence problems of stochastic mechanics and Bohmian mechanics have
two aspects in common with that of Bell's process: First, since the
law of motion (as defined by the drift in stochastic mechanics, the
velocity in Bohmian mechanics, and the jump rate in Bell's model) is
ill-defined at the nodes of the wave function, one needs to show that
the process never reaches a node. Second, while in stochastic
mechanics and Bohmian mechanics there are no jumps that could
accumulate, one needs to exclude (and has excluded) the analogous
possibility that the process could escape to infinity in finite time.

\section{The Result}
\label{sec:Result}

The basic ingredients of the model are:
\renewcommand{\labelitemi}{--}
\begin{itemize}
\item a complex Hilbert space $\Hilbert$ with inner product
$\sp{\,\cdot\,}{\,\cdot\,}$, the space of quantum states,
\item a self-adjoint operator $H$  acting on $\Hilbert$, the Hamiltonian,
\item an initial state vector $\Psi_0\in\Hilbert$  with $\|\Psi_0\|=1$,
\item a countable set $\conf$, physically thought of as configuration
   space and serving as state space of the jump process to be
   constructed, and
\item a positive-operator-valued measure (POVM) $\pov(\,\cdot\,)$ on $\conf$
acting on $\Hilbert$.
\end{itemize}
Here, a POVM is a family $(\pov(\x))_{x\in\conf}$ of positive bounded
self-adjoint operators on $\Hilbert$ such that, for each
$F\subseteq\conf$, the sum $\pov(F):=\sum_{x\in F}\pov(\x)$ exists in
the sense of the weak operator topology, and $\pov(\conf)=I$.  In
fact, the countable additivity then also holds in the strong topology
\cite{Davies}.  In particular,
\begin{equation}\label{L2sense}
    \forall\, \Phi \in \Hilbert: \: \sum_{\x \in \conf} \pov(\x) \Phi \text{
    converges in the $L^2$ sense to } \Phi.
\end{equation}
Every PVM is a POVM but not vice versa.  As has already been pointed
out in \cite{crea2}, the jump rate formula \eqref{tranrates} still
makes sense if $\pov(\,\cdot\,)$ is a POVM rather than a PVM.

In quantum field theory, the ``configuration observable''
$\pov(\,\cdot\,)$ is often a POVM; a typical situation is that
$\Hilbert$ is a subspace (e.g., the positive spectral subspace of the
free Hamiltonian) of a larger Hilbert space $\Hilbert_0$ containing
also unphysical states, and $\pov(\,\cdot\,) = P' \pov_0(\,\cdot\,)
I'$ where $P'$ is the projection $\Hilbert_0 \to \Hilbert$, $I'$ is
the embedding $\Hilbert \hookrightarrow \Hilbert_0$, and
$\pov_0(\,\cdot\,)$ is a PVM (the configuration observable) acting on
$\Hilbert_0$.

To establish the existence of a Markovian jump process with rates
\eqref{tranrates} we need the following joint assumption on $H$,
$\pov$, and the initial state vector $\Psi_0$.

\begin{ass}\label{ass}
The Hamiltonian $H$, the POVM $\pov$ and the state vector
$\Psi_0\in\Hilbert$ satisfy the conditions
\begin{enumerate}
\item[{\rm(A1)}] For all $t \in \RRR$ and $\x \in \conf$, $\Psi_t$ and
$ \pov(\x) \Psi_t$ belong to the domain of $H$.
\item[{\rm(A2)}] For all $ t_0, t_1 \in \RRR$ with $t_0 < t_1$,
\[
    \int_{t_0}^{t_1} dt \sum_{\x,\y \in \conf} \bigl| \sp{\Psi_t}
    {\pov(\y) H \pov(\x) \Psi_t} \bigr| < \infty \,.
\]
\end{enumerate}
\end{ass}

\bigskip

For given $H$ and $\pov$, Assumption \ref{ass} can also be understood
as an  assumption on $\Psi_0$, thus defining a set $\D \equiv \D_{H,\pov}
\subseteq \Hilbert$ of ``good'' state vectors for which the process is
well-defined. This $\D$ is invariant under the time evolution but not
necessarily a subspace of $\Hilbert$ because assumption (A2) is not
linear in $\Psi_t$. The following proposition provides conditions on
$H$ under which Assumption \ref{ass} holds for all $\Psi_0\in
\Hilbert$, so that $\D = \Hilbert$. (For general $H$ we do not know
how large $\D$ is, and whether  it is dense, as would be physically
desirable. We will not presuppose this but instead construct the
process solely for initial state vectors $\Psi_0\in\D$.)

\begin{proposition}\label{prop:D=H}
Assumption \ref{ass} holds for all $\Psi_0\in \Hilbert$ when
either
\begin{enumerate}
\item $H$ is bounded and $\conf$ is finite,  or
\item $H$ is a Hilbert--Schmidt operator, i.e.,  $\tr \, H^2 < \infty$.
\end{enumerate}
\end{proposition}

The proof is postponed until Section~\ref{sec: ProofProp}.
   Assumption (A1) implies that
$\pov(\y) H \pov(\x) \Psi_t$ exists, and thus that
\[
    \sigma_t(\y|\x) \text{ is well defined whenever }
    \sp{\Psi_t}{\pov(\x) \Psi_t} \neq 0.
\]
When $\sp{\Psi_t}{\pov(\x) \Psi_t} = 0$, we set $\sigma_t(\y|\x) :=
\infty$ for all  $\y$; thus, $\sigma_t(\y|\x)$ is always  defined as
a $[0,\infty]$-valued function.  (When suitably reinterpreted, the
numerator of \eqref{tranrates} still exists if $\pov(\x) \Psi_t$ and
$\pov(\y) \Psi_t$ merely lie in the  {form domain}, rather than
the domain, of $H$. We will not pursue here this kind of greater
generality.)

As was pointed out in the introduction, the rates $ \sigma_t(\y|\x)$ are
constructed in such a way that the corresponding Markov process $X_t$
should have the quantum distribution
\begin{equation}\label{equivar-distr}
\measure_t(\x) := \sp{\Psi_t}{\pov (\x) \Psi_t}\,,\quad \x\in\conf,
\end{equation}
at any time $t\in\R$. In other words, the family
$(\measure_t)_{t\in\R}$ should be \emph{equivariant}, or an
\emph{entrance law}, for the process.  Here is our main result stating
that such a process does exist.

\begin{theorem}\label{thm}
    Suppose Assumption \ref{ass} holds. Then there exists a
    right-continuous   (time-inhomogenous) Markovian pure jump process
    $(X_t)_{t\in\R}$ in $\conf$ with transition rates \eqref{tranrates} and
    such that, 
    for each $t$, $X_t$ has distribution $\measure_t$. The process
    is unique in distribution.
\end{theorem}

\section{The construction}

We fix some starting time $t_0\in\R$ and construct the process first
on the time interval $[t_0,\infty[$. We also introduce an auxiliary
``cemetery''  configuration $\friedhof$ in order to
deal with the possibility that the process explodes or runs into a
node. In the next section we will show that this does in fact not
occur. We write
\begin{equation}\label{node}
\node:= \{(t,\x)\in\R\times\conf: \measure_t(\x) = 0\}
\end{equation}
for the node-set of all exceptional times and positions for which the
transition rates \eqref{tranrates}  are infinite. Likewise,
\begin{equation}\label{conf_t}
\conf_t:= \{\x\in\conf: \measure_t(\x) > 0\}
\end{equation}
is the set of all admissible positions at time $t\in\R$.
Finally, for $(t,x)\notin\node$ we let
\begin{equation}\label{eq:tau}
    \tau_{t,x}:= \inf\{s>t: (s,x)\in\node\}
\end{equation}
be the first time instant after $t$ at which $x$ becomes a node; here
we set $\inf\emptyset:=\infty$.
Let us start with a technical lemma; its proof follows later in this section.
Formula \eqref{dPsi2dt} relies on our convention $\hbar=2$.

\begin{lemma}\label{measble}
   For every $\x\in \conf$,  the mapping $t \mapsto
   \measure_t(\x)$ is differentiable with locally integrable derivative
    \begin{equation}\label{dPsi2dt}
      \dot{\measure}_t(\x) = \Im \, \sp{\Psi_t}{\pov(\x) H \Psi_t}.
    \end{equation}
    In particular, the function $\measure_\cdot(\x)$ is locally
    absolutely continuous.  Also, the jump rates $\sigma_t(\y|\x)$
    depend measurably on $t$ with values in $[0,\infty]$, and the total
    jump rate
    \begin{equation}\label{sigmaconfdef}
      \g_\x(t) := \sum_{\y \in \conf} \sigma_t(\y|\x)
    \end{equation}
is finite whenever $(t,\x)\notin\node$.
\end{lemma}

The process $(X_t)_{t\ge t_0}$ will be constructed on the enlarged
position space $\conf\cup\{\friedhof\}$ by means of a suitable
sequence of random jump times $T_n$ and jump destinations $\dest_n$.
To achieve this we need two key quantities: the distribution
$\mu_{t,x}$ of the holding time in $\x\in\conf$
(i.e., the random waiting time before the next jump) after a given
time $t$, and the distribution $p_{t,x}$ of the the jump destination
at the jump time.  Our assumption that the process $(X_t)_{t\ge t_0}$
should have the transition rates \eqref{tranrates} simply means that
$\mu_{t,x}$ should be the distribution with ``failure rate function''
(or ``hazard rate function'') $\g_\x$; cf. e.g.  \cite[pp.~276 ff.,
577]{Ross}. That is, for any $(t,x)\notin\node$ we let $\mu_{t,x}$ be
the unique probability measure on $]t,\infty]$ with
``survival probabilities''
\begin{equation}\label{eq:mu}
    \mu_{t,\x}([u,\infty]) = \E^{ - \G_{t,\x}(u)} \text{ for all $u>t$,}
\end{equation}
where
\begin{equation}\label{eq:Gamma}
    \G_{t,\x}(u)=\int_{t}^u  \g_\x(s)\, ds\,;
\end{equation}
in \eqref{eq:mu} and below we set $\E^{-\infty}=0$.  (Note that
$\G_{t,\x}$ is left-continuous by the monotone convergence theorem, so
that there exists indeed a unique probability measure $\mu_{t,\x}$
having $\E^{ - \G_{t,\x}}$ as right-sided distribution function.)  In
particular, $\mu_{t,\x}(\{\infty\})>0$ if and only if
$\G_{t,\x}(\infty)<\infty$;  thus, in this case there is a
non-zero probability for the process to be frozen in $\x$.  If
$t=\infty$ or $\x=\friedhof$ we let $\mu_{t,\x}=\delta_\infty$ be the
Dirac measure at $+\infty$. The following lemma collects the essential
properties of $\mu_{t,\x}$.

\begin{lemma}\label{Lmudef}
    Suppose $(t,x)\notin\node$. Then the following statements hold:
    \begin{enumerate}
    \item On $]t,\tau_{t,x}[\,$, $\mu_{t,\x}$ has the density function
    $\g_\x \,\E^{ - \G_{t,\x}}$.
    \item If $\tau_{t,x}<\infty$ then $\mu_{t,\x}([\tau_{t,x},\infty])=0$.
    \item $\mu_{t,\x}(0<\g_\x<\infty)=1$.
    \item For any $u\in\,]t,\tau_{t,x}[\,$, 
    $\mu_{t,\x}(\, ]u,\infty]) >0$ and $\mu_{u,\x}=
    \mu_{t,\x}(\,\cdot\,|\, ]u,\infty])$.
    \end{enumerate}
\end{lemma}

The proof will be given later in this section.  It follows readily
from (a) that if we take $\mu_{t,\x}$ as the the distribution of the
holding time at $\x$ then the jump rate is indeed $\g_\x$, as
intended.  Assertion (b) states that the process cannot run into a
node by sitting on an $\x$ until it becomes a
node: right before $\x$ becomes a node, the rate $\g_\x$ grows so fast
that the process has probability 1 to jump away. In particular, the
unboundedness of the jump rates (even for bounded $H$, even for
Hilbert--Schmidt $H$) \emph{favours} the global existence rather than
preventing it. Statement (d) expresses a loss-of-memory property which
is responsible for the Markov property of the process.

As for the distribution of the jump destinations,
we obviously have to define
\begin{equation}\label{eq:destdistr}
p_{t,x}(\y) = \sigma_{t} (\y|\x)/\g_\x(t) \quad\text{ if }
0<\g_\x(t)<\infty\,.
\end{equation}
Otherwise we set $p_{t,x}=\delta_\friedhof$, the Dirac
measure at $\friedhof$.  The next lemma states that $p_{t,x}$ is
supported on the set $E_t$ of non-nodes defined in \eqref{conf_t}.

\begin{lemma}\label{LXdistribution}
$p_{t,x}(\conf_t)=1$ whenever 
$0<\g_\x(t)<\infty$.
\end{lemma}

With these ingredients we are now ready to construct the process
$(X_t)_{t\ge t_0}$ on $\conf\cup\{\friedhof\}$. Let
$(T_n,\dest_n)_{n\ge0}$ be a sequence of random variables with the
following properties. Let $T_0:=t_0$ and $\dest_0\in E_{t_0}$ be a
random variable with distribution $\measure_{t_0}$.  Then, for any
$n\ge0$, let $(T_{n+1},\dest_{n+1})$ have the conditional distribution
\begin{equation}\label{T+dest}
    \prob\big(T_{n+1} \in \,dt, \dest_{n+1}=\y\,| T_0,\dest_0,
\ldots, T_{n},\dest_{n}\big) = \mu_{T_n,\dest_n}(dt)\, p_{t,\dest_{n}}(\y)\,.
\end{equation}
The existence of such a sequence $(T_n,\dest_n)_{n\ge0}$ on a suitable
probability space $(\Omega_{t_0},\mathcal{F}_{t_0},\prob_{t_0})$
follows from the Ionescu--Tulcea theorem \cite[Theorem
5.17]{Kall}.  Moreover,
Lemmas~\ref{Lmudef} and~\ref{LXdistribution} imply that
$T_{n+1}<\tau_{T_n,Z_n}$ and $\dest_n\in\conf_{T_n}$ almost surely for
all $n$. We also have $T_n< T_{n+1}$ as long as $T_n< \infty$. So we
define
\begin{equation}\label{Xdef}
    X_t = \dest_n \text{ when }T_n\le t<T_{n+1}, \text{ and } X_t =\friedhof
     \text{ for }t\ge\zeta = \sup_n T_n\,.
\end{equation}
It is then clear that $(X_t)_{t\ge t_0}$ is right-continuous, and
$(t,X_t)\notin\node$ for all
$t\in[ t_0,\zeta[\,$ with probability $1$. In the
next section we will show that in fact $\zeta=\infty$ almost surely.

\bigskip
We now turn to the proofs of the lemmas above. Recall our convention that
$\hbar=2$.

\bigskip
\Proof[ of Lemma \ref{measble}] Since $\Psi_t$ belongs to the domain
of $H$ by Assumption (A1), the $\Hilbert$-valued mapping $t
\mapsto\Psi_t$ is differentiable with derivative $d\Psi_t/dt
=-\frac{\I}{2} H \Psi_t$ \cite[p.~265]{ReedSimon}. Hence
\[
    \big(\measure_{t+s}(\x) - \measure_t(\x)\big)/s = \sp{\Psi_{t+s}}
    {\pov(\x)(\Psi_{t+s} - \Psi_t)/s} + \sp{(\Psi_{t+s} - \Psi_t)/s}
    {\pov(\x) \Psi_t}
\]
converges, as $s \to 0$, to
\[
    \dot{\measure}_t(\x):=-\frac{\I}{2} \,\sp{\Psi_t} {\pov(\x)H \Psi_t}
    + \frac{\I}{2}\, \sp{H\Psi_t}{\pov(\x) \Psi_t} = \Im \,
    \sp{\Psi_t}{\pov(\x) H \Psi_t}\,.
\]
As a limit of the continuous difference ratios,
$t\mapsto\dot{\measure}_t(\x)$ is measurable. Moreover, using
\eqref{L2sense} we can write
\[
\begin{split}
\bigl|\dot{\measure}_t(\x)\bigr|&\le
\bigl|\sp{\Psi_t}{\pov(\x) H \Psi_t} \bigr|
= \bigl|\sp{H \pov(\x) \Psi_t}
    {\Psi_t} \bigr|\\
&= \Bigl|\sum_{\y \in \conf} \sp{H \pov(\x)\Psi_t}{\pov(\y) \Psi_t}
\Bigr| \leq \sum_{\y \in \conf} \bigl| \sp{\Psi_t}{\pov(\x) H \pov(\y)
\Psi_t} \bigr|\,.
\end{split}
\]
Together with Assumption (A2), it follows that $\dot{\measure}_t(\x)$
is locally integrable. In particular, $t \mapsto \measure_t(\x)$ is
locally absolutely continuous and  an integral function of $t
\mapsto\dot{\measure}_t(\x)$; see \cite[Theorem 6.3.10]{Cohn} or
\cite[Theorems 8.21 and 8.17]{rudin}.

Concerning the measurability of the jump rates $\sigma_t(\y|\x)$, it
is sufficient to show that $\sp{\Psi_t}{ \pov(\y) H \pov(\x) \Psi_t}$
depends measurably on $t$. (This is because $\measure_\cdot(\x)$ is
continuous and the ratio of nonnegative measurable functions is
measurable.)  To this end, we introduce the cutoff function $f_n(a):=
(a\wedge n)\vee(-n)$ and observe that, for every $t$,
\[
    \sp{\Psi_t}{ \pov(\y) f_n(H) \pov(\x) \Psi_t} \to \sp{\Psi_t}{
    \pov(\y) H \pov(\x) \Psi_t}
\]
as $n \to\infty$.  Hence $t \mapsto \sp{\Psi_t}{ \pov(\y) H \pov(\x)
\Psi_t}$ is a pointwise limit of continuous functions and thereby
measurable. In particular, the total jump rate $\g_\x$ is
measurable. For $(t,\x)\notin\node$ we have
\[
    \g_\x(t) \leq \sum_{\y\in\conf} \bigl| \sp{\Psi_t}{\pov(\y) H
    \pov(\x) \Psi_t} \bigr| \big/\measure_t(\x)\,.
\]
The last sum is finite because, due to \eqref{L2sense}, the series
$\sum _\y \sp{\Psi_t} {\pov(\y) H \pov(\x) \Psi_t}$ converges
to $\sp{\Psi_t} {H \pov(\x) \Psi_t}$ in every ordering,
and is therefore absolutely convergent.\EndProof

Before proving Lemma~\ref{Lmudef} we establish the following result, a
key fact for showing that the process never runs into a node. Recall
the definitions \eqref{eq:tau} and \eqref{eq:Gamma}.

\begin{lemma}\label{Lrateint}
    Suppose $(t,\x)\notin\node$. Then $\G_{t,\x}(u)<\infty$ if
    $t<u<\tau_{t,\x}$, while $\G_{t,\x}(\tau_{t,\x})=\infty$ if
    $\tau_{t,\x}<\infty$.
\end{lemma}

\Proof Consider first the case $t<u<\tau_{t,\x}$.  Since
$\measure_\cdot(\x)$ is continuous and positive on $[t,u]$, it stays
bounded away from zero on this interval.  On the other hand, we have
\[
    \int_{t}^{u} ds \sum_{\y\in \conf}\bigl[ \Im \, \sp{\Psi_s}{\pov(\y)
    H \pov(\x) \Psi_s}\bigr]^+ \leq \int_{t}^{u} ds \sum_{\y\in \conf}
    \bigl| \sp{\Psi_s}{\pov(\y) H \pov(\x) \Psi_s} \bigr| \,,
\]
and the last integral is finite due to Assumption (A2).
This proves the first assertion.

Consider now the case $\tau_{t,\x}<\infty$.
Since $\sum_i a_i^+ \geq\big[ \sum_i a_i \big]^+$ in general, we have
for all $t<s<\tau_{t,\x}$
\[
    \g_\x(s) \geq \big[ \Im \sum_{\y} \sp{\Psi_s}{\pov(\y) H \pov(\x)
    \Psi_s}\big]^+\big/ \measure_s(\x)\,.
\]
In view of \eqref{L2sense} and Lemma \ref{measble}, the last
expression is equal to
\[
\big[ \Im \,
    \sp{\Psi_s}{H \pov(\x) \Psi_s}\big]^+\big/\measure_s(\x)=
\big[- \dot{\measure}_s(\x)\big]^+  \big/\measure_s(\x)\,.
\]
Since always  $a^+ \geq a$, we arrive at the key inequality
\[
   \g_\x(s)  \geq  - \frac{d}{ds} \log \measure_s(\x)\,.
\]
The last derivative is integrable over any interval $[t,u]$ with
$t<u<\tau_{t,\x}$ because $\measure_s(\x)$ is bounded away from zero
on such an interval and $\dot{\measure}_s(\x)$ is locally integrable
by Lemma~\ref{measble}.  By the general fundamental theorem of
calculus as in \cite[Theorem 6.3.10]{Cohn} or \cite[Theorem
8.21]{rudin}, it follows that
\[
   \int_{t}^{u} \g_\x(s)\,ds \ge
   - \int_{t}^{u}  \frac{d}{ds} \log \measure_s(\x) \,ds
    =  \log \measure_t(\x) - \log \measure_u(\x).
\]
Letting $u\uparrow \tau_{t,\x}$ and using the continuity of $\measure_u(\x)$
we arrive at the second statement of the lemma.\EndProof

We are now ready for the proof of  Lemma \ref{Lmudef}.

\bigskip
\Proof[ of Lemma \ref{Lmudef}] (a) Let $t<u<\tau_{t,\x}$. Instead of
using the fundamental theorem of calculus (which would be possible),
we prefer to give here a direct argument which is based on Fubini's
theorem.  In view of Lemma \ref{Lrateint}, $\G_{t,\x}$ is finite on
$[t,u]$. Thus we can write, omitting the indices $t,\x$,
\[
\begin{split}
&\int_t^u \g(s)\,\E^{-\G(s)}ds = \int_t^u ds\,\g(s) \int_0^\infty
dr\,\E^{-r}\,1_{\{\G(s)\le r\}}\\
&=\int_0^\infty dr\,\E^{-r}\int_t^u  ds\,\g(s)\,1_{\{\G(s)\le r\}}
= \int_0^\infty dr\,\E^{-r}\ (r\wedge\G(u))\,.
\end{split}
\]
The last equality uses the fact that $\G$ is continuous and increasing.
Since $r\wedge\G(u)=r-[r-\G(u)]^+$, the last integral coincides with
$1- \E^{-\G(u)}=\mu_{t,\x}(]t,u])$, thus proving assertion (a).

(b) This is immediate from \eqref{eq:mu} and Lemma \ref{Lrateint}.

(c) This comes from statements (a) and (b) together with Lemma
\ref{measble}.

(d) Let $t<u<\tau_{t,x}$. Since $\G_{t,\x}(u)<\infty$ by Lemma \ref{Lrateint},
Equation \eqref{eq:mu} shows that  $\mu_{t,\x}(\, ]u,\infty]) >0$. Moreover,
for $v>u$ we have
\[
\mu_{t,\x}\big(\, ]v,\infty]\,\big|\, ]u,\infty]\big)=
\E^{-\G_{t,\x}(v)+\G_{t,\x}(u)}=\E^{-\G_{u,\x}(v)} =
\mu_{u,\x}\big(\, ]v,\infty]\,\big)
\]
by Equation \eqref{eq:Gamma}. This proves the final statement.\EndProof

We conclude this section with the proof of Lemma \ref{LXdistribution}.
\medskip

\Proof[ of Lemma \ref{LXdistribution}]
We only have to to show that $\sigma_t(\y|\x)=0$ whenever
$\measure_t(\y) =0$. But since $\|\pov(\y)^{1/2} \Psi_t\|^2=\measure_t(\y)$,
we then have $\pov(\y)^{1/2} \Psi_t=0$. Hence $\pov(\y)\Psi_t=0$
and therefore $\sp{\Psi_t}{\pov(\y) H \pov(\x) \Psi_t} = \sp{\pov(\y) \Psi_t}
{H \pov(\x) \Psi_t} =0$, which gives the result.\EndProof

\section{Non-explosion}

In the last section we have constructed a process $(X_t)_{t\ge t_0}$
that stays in the configuration space $\conf$ until some possibly
finite explosion  time $\zeta=\sup_n T_n$, at which it jumps
into the cemetery $\friedhof$.  We will now show that $\zeta$ is in
fact almost surely infinite. To this end we consider the random number
\begin{equation}\label{Sdef}
S(t) := \# \{n\ge 1: t_0 < T_n \leq t\} \in \ZZZ_+\cup\{\infty\}
\end{equation}
of jumps during the time interval $]t_0, t]$ for any $t>t_0$. We want
to show that $S(t)$ has finite expectation. To this end we start from
the following formula.

\begin{lemma}\label{LES} For all $t>t_0$,
\[
   \EEE_{t_0}\, S(t) = \int_{t_0}^{t} ds
    \sum_{\x,\y\in\conf}  \prob_{t_0}(X_s = \x)\,\sigma_s(\y|\x) \,.
\]
\end{lemma}

To estimate the last expression we will show:

\begin{lemma}\label{Lequi}
   $\prob_{t_0}(X_t =\x) \leq \measure_t(\x)$ for all $\x\in\conf$ and
   $t>t_0$.
\end{lemma}

In other words, though the rates are constructed in such a way that
the process should follow the equivariant distribution $\measure_t$,
we cannot exclude \emph{a priori} that some mass is lost at the
cemetery $\friedhof$.  Combining these two lemmas we obtain
\[
\begin{split}
    \EEE_{t_0}\, S(t) &\leq \int_{t_0}^{t} ds \sum_{\x,\y\in\conf}
    \measure_s(\x) \,  \sigma_s(\y|\x)\\
    & =  \int_{t_0}^{t} ds \sum_{\x,\y\in\conf}
    \bigl[\Im \, \sp{\Psi_s}{\pov(\y) H \pov(\x) \Psi_s} \bigr]^+  \\
    &  \leq  \int_{t_0}^{t} ds \sum_{\x,\y\in\conf}
    \bigl|\sp{\Psi_s}{\pov(\y) H \pov(\x) \Psi_s} \bigr|\,,
\end{split}
\]
and the last expression is finite by Assumption (A2).
Hence $S(t)<\infty$ almost surely, and thereby $\zeta>t$
almost surely. As $t$ was arbitrary, we conclude that $\zeta=\infty$
almost surely, as we wanted to show. We now turn to the proofs
of the two lemmas above.

\bigskip
\Proof[ of Lemma \ref{LES}] Using Equation \eqref{T+dest}
and Lemma \ref{Lmudef} we can write
\[
\begin{split}
\EEE_{t_0}\, S(t) &= \sum_{n\ge 0} \prob_{t_0}( t_0 \leq T_{n+1}\leq t)
   = \sum_{n\ge 0} \EEE \Bigl( \prob_{t_0}\big( t_0 \leq T_{n+1}
    \leq t\,|\, T_{k}, \dest_{k}:k\le n\big) \Bigr)
\\
    &= \sum_{n\ge 0} \ \EEE \  \int_{t_0}^{t} ds
    \, \1_{\{ T_n<s<\tau_{T_n,\dest_n} \}} \, \g_{\dest_{n}}(s) \, \E^{
    - \G_{T_n,\dest_{n}}(s)}
\\
    &=\sum_{n\ge 0} \int_{t_0}^{t} ds \; \EEE\Bigl(
    \1_{\{T_{n}< s\}} \, \g_{\dest_{n}}(s) \, \prob_{t_0}
    \bigl(T_{n+1}>s \,|\,  T_{k}, \dest_{k} :k\le n\bigr) \Bigr)
\\
   & = \int_{t_0}^{t} ds \sum_{\x\in \conf} \g_{\x}(s)
    \, \EEE \Bigl( \sum_{n\ge 0} \1_{\{T_{n}< t <T_{n+1},\,
\dest_{n}=\x\} } \Bigr)
\\
    &= \int_{t_0}^{t} ds \sum_{\x\in \conf} \g_{\x}(s)
    \, \prob_{t_0}( X_s=\x)\,.
\end{split}
\]
Together with \eqref{sigmaconfdef} the lemma follows.\EndProof

For the proof of Lemma \ref{Lequi} we consider the integral equation
\begin{equation}\label{int-eq}
    \rho_t(\x) = \measure_{t_0}(\x) \, \E^{-\G_{t_0,\x}(t)}
   + \sum_{\y\in \conf} \int_{t_0}^t ds \; \rho_s(\y)
    \, \sigma_s(\x|\y) \,  \E^{-\G_{s,\x}(t)}\,,
\end{equation}
$t\ge t_0,\,\x\in\conf$, for a time-dependent subprobability measure
$\rho_t$ on $\conf$.  Lemma \ref{Lequi} follows directly from the next
two results.

\begin{lemma}\label{lem: mini-solution}
The mapping $(t,\x)\mapsto \prob_{t_0}(X_t=\x)$ is the minimal
solution of \eqref{int-eq}.
\end{lemma}

\begin{lemma}\label{lem: equivar-solution}
The mapping $(t,\x)\mapsto \measure_t(\x)$ is a solution of
\eqref{int-eq} for arbitrary $t_0$.
\end{lemma}

\Proof[ of Lemma \ref{lem: mini-solution}]
For any $\x\in\conf$ and $t>t_0$
we can write
\[
   \prob_{t_0}(X_t = \x) = \sum_{n\ge 0} A_n(t,\x)\quad \text{ with } \
A_n(t,\x):=\prob_{t_0}\bigl( T_n \leq t < T_{n+1},\, \dest_n =\x\bigr).
\]
It follows from \eqref{eq:mu} that
\begin{equation}\label{A0expression}
    A_0(t,\x)=\measure_{t_0}(\x)  \, \E^{-\G_{t_0,\x}(t)}
\end{equation}
and, for $n\ge1$,
\[
\begin{split}
   A_n(t,\x)&=
      \sum_{\x_0, \ldots, \x_{n-1} \in \conf} \quad
       \idotsint\limits_{t_0<t_1<\cdots <t_n\le t}\\
      &\quad \times \prob_{t_0} \Bigl( T_1 \in dt_1, \ldots,
      T_n \in dt_n, T_{n+1} > t,
      \dest_0 = \x_0, \ldots, \dest_n =\x_n \Bigr)
\\
      &=
      \sum_{\x_0, \ldots, \x_{n-1} \in \conf} \quad
      \idotsint\limits_{t_0<t_1<\cdots <t_n\le t}dt_1\cdots dt_n
\\
      &\quad\times \measure_{t_0}(\x_0)\, \Bigl(\,
      \prod_{i=1}^{n} \E^{-\G_{t_{i-1},\x_{i-1}}(t_{i})}
\sigma_{t_i}(\x_i|\x_{i-1}) \,
      \Bigr) \,   \E^{-\G_{t_n,\x}(t)} \,,
\end{split}
\]
where $\x_n:=\x$. In particular, separating the summation over $\x_{n-1}$ and
the integration over $t_n$ we find that
\begin{equation}\label{Anexpression}
    A_n(t,\x)=\sum_{y\in\conf} \int_{t_0}^t ds\ A_{n-1}(s,\y)\,
    \sigma_s(\x|\y) \, \E^{-\G_{s,\x}(t)}\,.
\end{equation}
This shows that $ \prob_{t_0}(X_t = \x)$ satisfies \eqref{int-eq}.

Now let $\rho_t(\x)$ be an arbitrary (nonnegative) solution of
\eqref{int-eq}.  An $(N-1)$-fold iteration of \eqref{int-eq} then
leads to the equation
\[
   \rho_t(\x) = \sum_{n=0}^{N-1} A_n(t,\x) + R_N(t,\x),
\]
with $A_n(t,\x)$ defined by \eqref{A0expression}
and \eqref{Anexpression}, and the remainder term
\[
\begin{split}
R_N(t,\x)&= \sum_{\x_0, \ldots, \x_{N-1} \in \conf}\quad
     \idotsint\limits_{t_0<t_1<\cdots <t_N\le t}dt_1\cdots dt_N
\\
    &\qquad \times  \rho_{t_1}(\x_0) \, \Bigl( \prod_{i=1}^{N-1}
    \sigma_{t_i}(\x_i|\x_{i-1}) \,
     \E^{-\G_{t_{i},\x_{i}}(t_{i+1})} \Bigr) \,
   \sigma_{t_{N}}(\x|\x_{N-1}) \,  \E^{-\G_{t_N,\x}(t)}\,.
\end{split}
\]
(Compared with $A_N(t,\x)$, $R_N(t,\x)$ involves $\rho_{t_1}$ rather
than $\rho_{t_0}=\measure_{t_0}$, and the $\sigma$'s and $\E^{-\G}$'s
run in a different order.)  Since $R_N(t,\x)\ge0$, we see that $
\rho_t(\x)$ exceeds each partial sum of the infinite series
constituting $\prob_{t_0}(X_t=\x)$. This proves Lemma \ref{lem:
  mini-solution}.\EndProof

\Proof[ of Lemma \ref{lem: equivar-solution}] We start from the
observation that, by \eqref{dPsi2dt}, \eqref{L2sense},
\eqref{tranrates} and the self-adjointness of $H$ and $\pov(\x)$,
\[
\begin{split}
    \dot{\measure}_t(\x) & =\Im \, \sp{H \pov(\x) \Psi_t}{\Psi_t}
= \sum_{\y\in \conf}
      \Im \, \sp{H \pov(\x) \Psi_t}{\pov(\y) \Psi_t}
\\
    &= \sum_{\y \in \conf} \Bigl( \measure_t(\y)\, \sigma_t(\x|\y)  -
    \measure_t(\x)\, \sigma_t(\y|\x)  \Bigr)
\\
    &= \sum_{\y \in \conf}
   \measure_t(\y)\, \sigma_t(\x|\y) - \measure_t(\x)\, \g_\x(t) \,
\end{split}
\]
This means that the integral equation \eqref{int-eq} for $\measure_t(\x)$
takes the form
\begin{equation}\label{eq:part-integr}
    \measure_{t}(\x)  - \measure_{t_0}(\x) \, \E^{-\G_{t_0,\x}(t)}
   = \int_{t_0}^{t} ds \,\bigl(
      \dot{\measure}_s (\x) + \measure_s(\x) \,
      \g_\x(s)  \bigr)  \E^{-\G_{s,\x}(t)}.
\end{equation}
To establish this equation we  write for brevity
$f(s)=\measure_s (\x)$ and $g(s)=\E^{-\G_{s,\x}(t)}$ and distinguish
two cases.

\emph{Case 1:} $f>0$ on $[t_0,t]$; that is, $\x$ is never a node on
this interval.  Then, by Lemma~\ref{Lrateint}, $\g_\x$ is integrable
over $[{t_0},t]$, whence $s\mapsto\G_{s,\x}(t)$ is absolutely
continuous with derivative $-\g_\x$.  Since the exponential function
is Lipschitz on $]{-}\infty,0]$, it follows that $g$ is absolutely
continuous with derivative $\dot{g}=\g_\x\, g$
Lebesgue-almost-everywhere; see \cite[Corollary 6.3.7]{Cohn} or
\cite[Theorem 8.17]{rudin}.  Equation \eqref{eq:part-integr} is thus
equivalent to the partial integration formula
\[
    f(t)g(t)- f(t_0)g(t_0) = \int_{t_0}^t\,\bigl( \dot{f}(s)g(s) +
    f(s)\dot{g}(s)\bigr) ds
\]
which holds according to Corollary 6.3.8 of \cite{Cohn}.

\emph{Case 2:} $f(s)=0$ for some $s\in[t_0,t]$; that is, $\x$ is a
node at some time $s$. By the continuity of $f$, there exists then a
largest such $s$ in $[t_0,t]$, say $\theta$. Suppose first that
$\theta=t_0$.  We can then apply Case 1 to each subinterval $[t_*,t]$
of $[t_0,t]$, which yields \eqref{eq:part-integr} with $t_*$ in place
of $t_0$.  Let $I(t_*)$ be the corresponding integral on the
right-hand side. Since the integrand is nonnegative, we can use the
monotone convergence theorem to conclude that $I(t_*)\uparrow I(t_0)$
as $t_*\downarrow t_0$. On the other hand, $f(t_*)g(t_*)\to
0=f(t_0)g(t_0)$ as $t_*\downarrow t_0$ because $f$ is continuous and
$0\le g\le1$. This proves \eqref{eq:part-integr} in the case
$\theta=t_0$.

If $\theta>t_0$, we observe that $\G_{s,\x}(t)=\infty$ for all
$s<\theta$. Indeed, we even have $\G_{s,\x}(\theta)= \infty$. This is
evident when $]s,\theta[$ consists only of node-times because then
$\g_x$ is infinite on this interval; otherwise it follows from the
second statement of Lemma~\ref{Lrateint} applied to the segment from a
non-node-time between $t$ and $\theta$ to the next node-time.
Consequently, the left-hand side of \eqref{eq:part-integr} is equal to
$f(t)$, while the integrand on the right-hand side vanishes on
$[t_0,\theta[$.  This means that we have to establish
\eqref{eq:part-integr} with $t_0$ replaced by $\theta$.  But this is
trivial when $\theta=t$ because then both sides vanish, and otherwise
follows from the previous paragraph.\EndProof

It is now easy to complete the proof of the theorem.

\medskip \Proof[ of Theorem \ref{thm}] As we have shown above, for any
$t_0\in\R$ there exists a right-continuous pure jump process
$(X_t)_{t\ge t_0}$ on a suitable probability space
$(\Omega_{t_0},\mathcal{F}_{t_0},\prob_{t_0})$.  Since $\zeta=\infty$
almost surely, this process avoids the cemetery $\friedhof$ and thus
takes values in $\conf$. Hence $\sum_{\x \in \conf} \prob_{t_0}(X_t
=\x ) =1$ for all $t\ge t_0$.  Lemma \ref{Lequi} therefore implies
that $\prob_{t_0}(X_t =\x ) = \measure_t(\x)$ for all $\x\in\conf$ and
$t>t_0$.  In particular, if $\conf_t$ is given by \eqref{conf_t} then
$X_t\in\conf_t$ for all $t\ge t_0$ with probability 1.

We also note  that $(X_t)_{t\ge t_0}$ is Markovian; its transition
matrix from time $s$ to time $t$ given by
\[
P_{s,t}(\x,\cdot)=\prob_s(X_t=\cdot\,|\,X_s=\x)
\]
when $\x\in\conf_s$, and arbitrary otherwise. This follows directly
from the construction together with Lemma \ref{Lmudef}(d). In
particular, the distribution $\rmP_{t_0}$ of $(X_t)_{t\ge t_0}$ on the
Skorohod space $D([t_0,\infty[,\conf)$ of all
c\`adl\`ag\footnote{continues \`a droite avec des limites \`a gauche =
  right-continuous with left limits.}  paths from $[t_0,\infty[$ to
$\conf$ is uniquely determined, and the family
$(\rmP_{t_0})_{t_0\in\R}$ is consistent. Kolmogorov's extension
theorem \cite[Theorem 5.16]{Kall} therefore provides us with a
probability measure $\rmP$ on $E^\R$ which extends all distributions
$\rmP_{t_0}$ and is therefore concentrated on $D(\R,\conf)$, the space
of all c\`adl\`ag paths on $\R$.  Under $\rmP$, the canonical
coordinate process constitutes the global Markov jump process with the
desired properties.\EndProof

\section{Proof of Proposition \ref{prop:D=H}}
\label{sec: ProofProp}

First we consider case (a). Since $H$ is bounded, assumption (A1)
holds trivially. The boundedness of $H$ also implies that the
expression $\sp{\Psi_t} {\pov(\x) H \pov(\y) \Psi_t}$ is (well defined
and) a continuous function of $t$ for every $\x$ and $\y$.
As $\conf$ is finite, the integrand in
Assumption (A2) is continuous and therefore locally integrable.

Turning to case (b), we observe first that assumption (A1)
is again trivially satisfied because
Hilbert--Schmidt operators are bounded.
Assumption (A2) will follow from the inequality
\begin{equation}\label{HSineq}
    \sum_{\x,\y \in \conf} \bigl| \sp{\Psi} {\pov(\x) H
    \pov(\y) \Psi} \bigr| \leq \|\Psi\|^2 \sqrt{\tr \, H^2} \quad
    \forall\, \Psi \in \Hilbert
\end{equation}
which we prove now.

We start with a general remark. Let $\Indexset$ be a countable index
set and $A_i$ and $B_i$, $i\in \Indexset$, any Hilbert--Schmidt
operators with (possibly different) adjoints $A_i^*$ resp.\ 
$B_i^*$; i.e., we have $\tr \, A_i^* A_i < \infty$ and similarly for
$B_i$.  The Cauchy--Schwarz inequality then asserts that
\begin{equation}\label{strongCSineq}
    \sum_{i\in \Indexset} |\tr \, A_i^* B_i| \leq
    \Bigl(\sum_{i\in \Indexset} \tr \, A_i^* A_i \Bigr)^{1/2} \,
    \Bigl(\sum_{i\in \Indexset} \tr \, B_i^* B_i \Bigr)^{1/2}
\end{equation}
whenever both terms on the right hand side are finite. (Note that we can
put the modulus sign inside of the sum because we can replace
$A_i$ by $z_i A_i$ with $z_i = (\tr \, A_i^* B_i)/|\tr \, A_i^* B_i|$
whenever $\tr\, A_i^* B_i \neq 0$.)

To obtain \eqref{HSineq} from \eqref{strongCSineq}, we set $\Indexset
= \conf \times \conf$, $A_{\x,\y} = \pov(\x)^{1/2} P_\Psi
\pov(\y)^{1/2}$ with $P_\Psi = |\Psi \rangle
\langle \Psi |$ the projection to $\CCC\Psi$, and $B_{\x,\y} =
\pov(\x)^{1/2} H \pov(\y)^{1/2}$. Then $ \tr \, A^*_{\x,\y} B_{\x,\y}
=\sp{\Psi}{\pov(\x) H \pov(\y) \Psi}.$ To see that $A_{\x,\y}$ is a
Hilbert--Schmidt operator, we note that
\[
    \tr \, A_{\x,\y}^* A_{\x,\y} = \tr \Bigl( \pov(\y)^{1/2} P_\Psi
    \pov(\x) P_\Psi \pov(\y)^{1/2} \Bigr) = \sp{\Psi}{\pov(\x) \Psi }
    \sp{\Psi}{\pov(\y) \Psi} < \infty.
\]
It follows further from \eqref{L2sense} that
\[
    \sum_{\x,\y} \tr \, A^*_{\x,\y} A_{\x,\y} =\sum_{\x,\y} \sp{\Psi}{
    \pov(\x) \Psi } \sp{\Psi}{\pov(\y) \Psi}= \|\Psi\|^4.
\]

Next we show that $B_{\x,\y}$ is a Hilbert--Schmidt operator. Note
that $0\leq \pov(\x) \leq I$  since $I - \pov(\x) = \pov(\conf
\setminus \{\x\}) \geq 0$. This implies that $0 \leq
\sp{\Phi}{\pov(\x) \Phi} \leq \sp{\Phi}{\Phi}$ for all $\Phi \in
\Hilbert$. Setting $\Phi := C\phi_n$ for any Hilbert--Schmidt operator
$C$ and an orthonormal basis $\{\phi_n: n \in \NNN\}$ of $\Hilbert$ we
find
\[
    \sum_{n \in \NNN} \sp{\phi_n}{C^* \pov(\x) C \phi_n} \leq \sum_{n \in
\NNN}
     \sp{\phi_n}{C^* C \phi_n}
\]
and thus
\begin{equation}\label{CpovC}
    \tr \, C^* \pov(\x) C \leq \tr \, C^* C.
\end{equation}
That is, if $C$ is a Hilbert--Schmidt operator then so is
$\pov(\x)^{1/2} C$; and so is $C \pov(\x)^{1/2} = (\pov(\x)^{1/2}
C^*)^*$. As a consequence, $B_{\x,\y}$ is a Hilbert--Schmidt operator.

Finally, we need to show that
\[
    \sum_{\x,\y \in \conf} \tr \, B^*_{\x,\y} B_{\x,\y} \leq \tr \, H^2.
\]
For every finite subset $F \subseteq \conf$ we have, using the linearity
of the trace and its invariance under cyclic permutations,
\[
\begin{split}
    &\sum_{\x,\y \in F} \tr \, B^*_{\x,\y} B_{\x,\y} =
    \sum_{\x,\y \in F} \tr \, H\pov(\x) H \pov(\y)
    = \tr \, H \pov(F) H \pov(F)\\
    &= \tr \bigl(
    \pov(F)^{1/2} H \pov(F)^{1/2} \bigr)^* \bigl(
    \pov(F)^{1/2} H \pov(F)^{1/2} \bigr)
    \leq \tr \, H^2\,.
\end{split}
\]
The last inequality comes from the fact that, according to
\eqref{CpovC}, the Hilbert--Schmidt norm $\|C\|_{HS} = \sqrt{\tr \,
C^* C}$ of an operator $C$ can only decrease when $C$ is multiplied,
from the left or from the right, by $P^{1/2}$ where $0 \leq P \leq
I$. Taking the supremum over all finite subsets
$F$ and combining all inequalities above we arrive at \eqref{HSineq}.

\bigskip
\noindent \textit{Acknowledgements.}  We thank Guido Bacciagaluppi of
the University of California, Berkeley, for his friendly
correspondence and Sheldon Goldstein of Rutgers University for helpful
discussions.  R.T.\ gratefully acknowledges support by the German
National Science Foundation (DFG), and hospitality at the Mathematics
Department of Rutgers University.

\end{document}